\documentclass[review,12pt]{elsarticle}

\usepackage[a4paper,width=170mm,top=30mm,bottom=25mm]{geometry}
\usepackage{lineno,hyperref}
\usepackage{multirow}
\usepackage{amsmath}
\numberwithin{figure}{section}
\numberwithin{equation}{section}
\numberwithin{table}{section}
\usepackage{color}
\usepackage{booktabs}
\usepackage{graphicx}
\usepackage{subfigure}
\usepackage{setspace}
\usepackage{amssymb}
\usepackage{bm}

\journal{arXiv Preprint}

\begin{document}

\begin{frontmatter}

\title{Two-dimensional DtN-FEM scattering analysis of SH guided waves by an interface debonding in a double-layered plate}


\author[mymainaddress,mysecondaryaddress,mythirdaddress]{Chen Yang\corref{mycorrespondingauthor}}
\ead{yangc@pku.edu.cn; yc@nuaa.edu.cn; yang.c.ag@m.titech.ac.jp}

\author[mysecondaryaddress]{Ruigang Qin\corref{mycorrespondingauthor}}
\ead{q827135409@outlook.com}

\author[mysecondaryaddress]{Sohichi Hirose}
\ead{shirose@cv.titech.ac.jp}

\author[mythirdaddress]{Bin Wang}
\ead{wangbin1982@nuaa.edu.cn}

\author[mythirdaddress]{Zhenghua Qian}
\ead{qianzh@nuaa.edu.cn}

\cortext[mycorrespondingauthor]{Corresponding author}
\address[mymainaddress]{State Key Lab for Turbulence and Complex Systems, College of Engineering, Peking University, Beijing 100871, China}
\address[mysecondaryaddress]{Tokyo Institute of Technology, 2-12-1-W8-22, O-okayama, Meguro, Tokyo, 152-8552, Japan}
\address[mythirdaddress]{Nanjing University of Aeronautics and Astronautics, 29 Yudao Jie, Nanjing 210016, China}

\begin{abstract}

In this paper, a two-dimensional Dirichlet-to-Neumann (DtN) finite element method (FEM) is developed to analyze the scattering of SH guided waves due to an interface delamination in a bi-material plate. During the finite element analysis, it is necessary to determine the far-field DtN conditions at virtual boundaries where both displacements and tractions are unknown. In this study, firstly, the scattered waves at the virtual boundaries are represented by a superposition of guided waves with unknown scattered coefficients. Secondly, utilizing the mode orthogonality, the unknown tractions at virtual boundaries are expressed in terms of the unknown scattered displacements at virtual boundaries via scattered coefficients. Thirdly, this relationship at virtual boundaries can be finally assembled into the global DtN-FEM matrix to solve the problem. This method is simple and elegant, which has advantages on dimension reduction and needs no absorption medium or perfectly matched layer to suppress the reflected waves compared to traditional FEM. Furthermore, the reflection and transmission coefficients of each guided mode can be directly obtained without post-processing. This proposed DtN-FEM will be compared with boundary element method (BEM), and finally validated for several benchmark problems.

\end{abstract}

\begin{keyword}
SH guided waves; Interface delamination; DtN-FEM; Virtual boundaries; Mode orthogonality
\end{keyword}

\end{frontmatter}


\section{Introduction}

Nondestructive testings demonstrate its absolute advantages and potential compared with other traditional destructive methods over the past few decades. As one of the nondestructive testings, the ultrasonic guided waves have naturally captured the attention of researchers. Although the conventional bulk wave detection technique has been widely used in many fields, the guided wave has lower preconditions, is more convenient to operate and needs to process fewer waveforms. These unique superiors cannot be replicated by bulk wave detection. Undoubtedly, the emergence of the guided wave technique has considerably reduced the cost of nondestructive testing.

Currently, many researchers have obtained huge achievements in quantitative non-destructive testing (QNDT) using ultrasonic guided waves \cite{Da20201479,Da20201923}. Wang proposed a prevalent method for shape and depth construction of plate, half space and layered semi-infinite space using whole reflection coefficients of SH waves \cite{Wang20121782,Wang2015}, Rayleigh Wave \cite{Wang2018} and Love waves \cite{Wang2019}, respectively. In addition, Da et al. \cite{Da2018181} utilize scattering coefficients of the torsional modes to reconstruct the defect shape in pipeline structures effectively. Consequently, it is enormously important to require a thorough understanding and accurate computation of forward scattering phenomenon, to acquire a database for QNDT research such as near-and-far field scattering data.

Graciously according to the polarizing direction, guided waves propagating in a plate can be roughly divided into two categories, which in common are Lamb wave and SH wave. This study uses an SH guided wave to accurately detect the interface delamination. SH wave which is an out-plane polarization guided wave, not purely has simple mechanical properties, but also has obvious mode conversion, and less energy loss when encountering obstacles. Because of these fundamental mechanical properties, it is suitable to detect damage in many fields.

In the analysis of elastic guided wave scattering problems, two key aspects must be properly addressed and treated. First, due to the dispersive and multi-mode nature of guided and evanescent waves, the wavenumber and wave structures should be accurately determined by dispersion equations. Moreover, dispersion curves can be drawn analytically by solving transcendent equations with complex root-racking modulus-converging algorithm \cite{Yang2018}. This information is useful for mode selection, as well as generation and reception of a single mode in non-destructive evaluation. Second, SH wave mode interaction with defects must be accurately determined, which is the focus of this paper. The complexity of mode conversion restrains an analytical approach to very simple geometries, and for more general cases, it requires valid numerical methods like FEM and BEM to obtain an accurate scattering wave-field. In our previous work, a modified BEM was proposed to solve the three-dimensional plate scattering problem \cite{Yang2021145}, where only the interfaces and flaw boundaries need to be discretized. However, BEM needs to handle the singularity and compute the fundamental solution matrix numerically, and the final BEM matrix will be full-rank, which consume many storage and computing resources, especially when simulating high frequency scattering problems. 

Compared to BEM, FEM is particularly attractive for forward scattering analysis since FEM is more simple and popular, and the final global FEM matrix is sparse. For plate scattering computation, Koshiba \cite{Koshiba198418} implemented the analysis of the scattering problem of Lamb waves in an elastic waveguide by FEM. Al-Nassar and Datta \cite{Al-Nassar1991125} used the Lamb wave to detect a normal rectangular strip weldment and obtained the guided waves’ behavior at specific frequencies. A semi-analytical FEM proposed by Hayashi and Rose \cite{Hayashi200375} was used to analyze the behavior of guided waves in plates and pipes. This method was soon applied to the field of flaw detection, for example, rail base \cite{Ramatlo2020,Yue2020}. Gunawan and Hirose \cite{Gunawan2004996} provide a mode-exciting method to study the scattering phenomenon in an infinite plate. Moreover, the researchers also conducted related explorations in the processing approach of the scattering boundary. The orthogonality of guided wave modes was provided to solve the scattering problem in a plate \cite{Shkerdin20042089}. Other researchers, such as Moreau \cite{Moreau2006611}, and Kubrusly \cite{Kubrusly2019,Kubrusly2021} try to apply this method to study relative objects.

The goal of this paper is to propose a FEM implemented in a bi-material plate to avoid the computation error caused by spurious reflected waves in 2-D scattering problems. During the finite element analysis, it is necessary to determine the far-field DtN conditions at virtual boundaries where both displacements and tractions are unknown. In this study, firstly, the scattered waves at the virtual boundaries are represented by a superposition of guided waves with unknown scattered coefficients. Secondly, utilizing the mode orthogonality, the unknown tractions at virtual boundaries are expressed in terms of the unknown scattered displacements at virtual boundaries via scattered coefficients. Thirdly, this relationship at virtual boundaries can be finally assembled into the global DtN-FEM matrix to solve the problem. This method is simple and elegant, which has advantages on dimension reduction and needs no absorption medium or perfectly matched layer to suppress the reflected waves compared to traditional FEM. Furthermore, the reflection and transmission coefficients of each guided mode can be directly obtained without post-processing. This proposed DtN-FEM will be compared with boundary element method (BEM), and validated for several benchmark problems. Finally, this proposed DtN-FEM will be applied into investigating the effect on different materials, various lengths of the delamination and location of the interface.

\section{DtN FEM}
\subsection{General description of the problem}
\begin{figure}[h]
	\centering
	\includegraphics[scale=0.8]{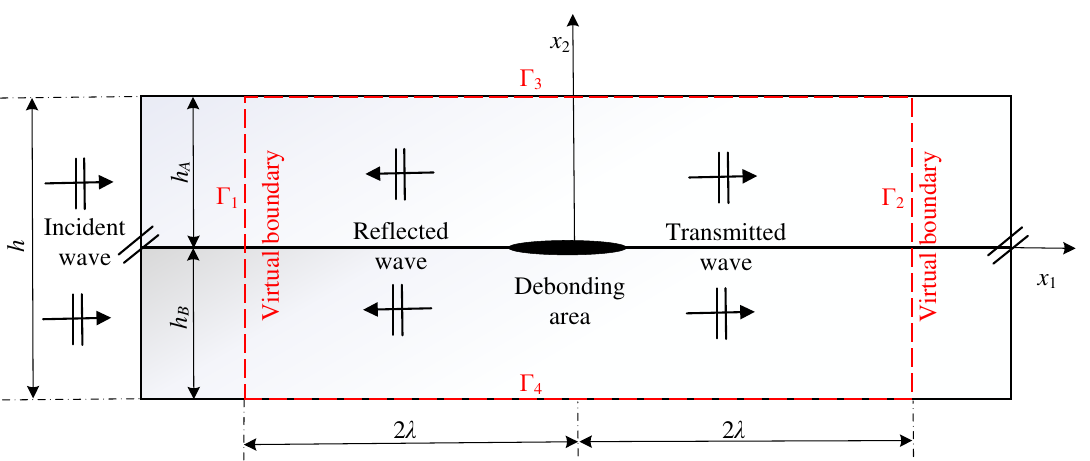}
	\caption{An infinite double-layered isotropic plate is considered where the material of upper layer $A$ is steel and the material of lower layer $B$ can be steel/aluminum/titanium, $h_A+h_B=h$, the distance between virtual boundaries and debonding area is at least $2\lambda$ where $\lambda$ is the wavelength of the first SH guided mode, debonding area is located at the interface and an incident SH wave is excited in far-field, traveling as a guided wave, and scattered by an interface debonding.}
	\label{fig_problem}
 \end{figure}

\begin{figure}[h]
	\centering
	\includegraphics[scale=0.8]{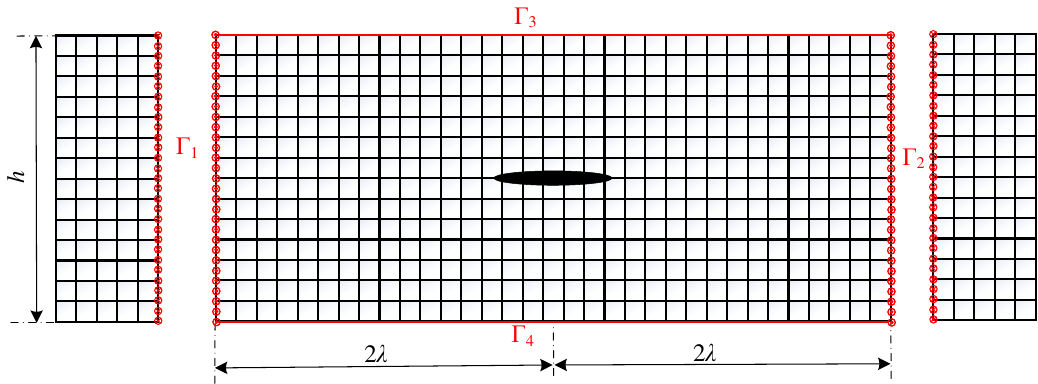}
	\caption{The domain mesh where displacements and tractions at the virtual boundaries $\Gamma _1$ and $\Gamma _2$ are both unknown.}
	\label{fig_mesh}
 \end{figure}

An infinite double-layered isotropic plate is considered as shown in Fig.\ref{fig_problem}, where the material of upper layer $A$ is steel and the material of lower layer $B$ can be steel/aluminum/titanium. The thicknesses of layer $A$ and layer $B$ are $h_A$ and $h_B$ respectively, and debonding area is located at the interface. Plane stress deformation in the plane of $x_1-x_2$ will be considered here. The upper and lower surfaces of the plate and the debonding surface are stress-free. Let us consider an incident SH wave excited in far-field, traveling as a guided wave, and scattered by a debonding surface. The aim of this paper is to simulate the scattered field. In traditional FEM, extra absorbing materials or perfectly matched layered is utilized to absorb the spurious reflected waves in order to simulate the infinite domain. However, in this paper, the domain surrounded by the red dashed line in Fig.\ref{fig_problem} is the actual simulation area, and the domain mesh is demonstrated in Fig.\ref{fig_mesh}. The key to the problem is how to deal with the far-field DtN conditions at the virtual boundaries since displacements and tractions at the virtual boundaries are both unknown.

\subsection{The mode orthogonality}

In order to get the relationship between displacements and tractions at the virtual boundary, the mode orthogonality should be proved firstly. Now consider two states, which are two guided SH modes with different wavenumber $k_m$ and $k_n$, and incident in the positive direction of axis $x_1$ in a double-layered isotropic plate without flaws. In accordance with the displacement patterns and stress patterns of 2D guided waves, we describe incident waves as

\begin{equation}
u_i^m\left( {{x_1},{x_2}} \right) = \hat u_i^m\left( {{x_2}} \right){\text{e}^{{\rm{i}}{k_m}{x_1}}},\quad \sigma _{ij}^m\left( {{x_1},{x_2}} \right) = \hat \sigma _{ij}^m\left( {{x_2}} \right){\text{e}^{{\rm{i}}{k_m}{x_1}}}.
\label{eq2.1}
\end{equation}
\noindent 
and
\begin{equation}
u_i^n\left( {{x_1},{x_2}} \right) = \hat u_i^n\left( {{x_2}} \right){\text{e}^{{\rm{i}}{k_n}{x_1}}},\quad \sigma _{ij}^n\left( {{x_1},{x_2}} \right) = \hat \sigma _{ij}^n\left( {{x_2}} \right){\text{e}^{{\rm{i}}{k_n}{x_1}}}.
\label{eq2.2}
\end{equation}

\noindent
which can be understood as expanding both displacement form and stress form as the product of the amplitude term and the propagating term.

Applying reciprocal theorem and considering the closed red line as the integral boundary shown in Fig., the following BIE can be obtained:

\begin{equation}
\int\limits_{{\Gamma _{1 - 4}}} {\hat u_i^m\left( {{x_2}} \right){\text{e}^{{\rm{i}}{k_m}{x_1}}}\hat \sigma _{ij}^n\left( {{x_2}} \right){\text{e}^{{\rm{i}}{k_n}{x_1}}}{n_j}} \text{d}\Gamma \left( \bm{x} \right)\; - \int\limits_{{\Gamma _{1 - 4}}} {\hat u_i^n\left( {{x_2}} \right){\text{e}^{{\rm{i}}{k_n}{x_1}}}\hat \sigma _{ij}^m\left( {{x_2}} \right){\text{e}^{{\rm{i}}{k_m}{x_1}}}{n_j}} \text{d}\Gamma \left( \bm{x} \right)\; = 0
\label{eq2.3}
\end{equation}

\noindent
where $n_k$ is the \textit{k}th component of unit vector outward surface normal to $\Gamma $.

Given the traction free boundary conditions on the upper and lower surfaces, the above equation can be simplified as

\begin{equation}
\left( {{\text{e}^{{\rm{i}}{k_m}{x_b}}}{\text{e}^{{\rm{i}}{k_n}{x_b}}} - {\text{e}^{{\rm{i}}{k_m}{x_a}}}{\text{e}^{{\rm{i}}{k_n}{x_a}}}} \right) \cdot {\int\limits_{-h_B}^{{h_A}} {\left[ {\hat u_3^m\left( {{x_2}} \right)\hat \sigma _{31}^n\left( {{x_2}} \right) - \hat u_3^n\left( {{x_2}} \right)\hat \sigma _{31}^m\left( {{x_2}} \right)} \right]\text{d}{x_2}} }  = 0
\label{eq2.4}
\end{equation}

It is clearly found that when $m=n$, the product of the second term in Eq.\eqref{eq2.4} automatically satisfies and when $m\ne n$, we can get

\begin{equation}
{\int\limits_{-h_B}^{{h_A}} {\left[ {\hat u_3^m\left( {{x_2}} \right)\hat \sigma _{31}^n\left( {{x_2}} \right) - \hat u_3^n\left( {{x_2}} \right)\hat \sigma _{31}^m\left( {{x_2}} \right)} \right]\text{d}{x_2}} }  \equiv 0
\label{eq2.5}
\end{equation}

\noindent
since the variables $x_a$ and $x_b$ are both arbitrary.

Because $\sigma _{31}^{n,m}\left( {{x_2}} \right)$ can be expressed as

\begin{equation}
\sigma _{31}^{n,m}\left( {{x_2}} \right) = \mu {\rm{i}}{k_{n,m}}\hat u_3^{n,m}\left( {{x_2}} \right)
\label{eq2.6}
\end{equation}

Eq.\eqref{eq2.5} can be simplified as

\begin{equation}
\int\limits_{-h_B}^{{h_A}} { {\hat u_3^m\left( {{x_2}} \right)\hat u_3^n\left( {{x_2}} \right)\text{d}{x_2}} }  \equiv 0
\label{eq2.7}
\end{equation}

The same conclusion can be drawn when we consider two states which are two guided modes with different wavenumbers, and incident in different directions. Considering an incident wave of wavenumber $k_m$ propagating in the negative direction of axis $x_1$

\begin{equation}
u_i^m\left( {{x_1},{x_2}} \right) = \hat u_i^{ - m}\left( {{x_2}} \right){\text{e}^{ - {\rm{i}}{k_m}{x_1}}},\quad \sigma _{ij}^m\left( {{x_1},{x_2}} \right) = \hat \sigma _{ij}^{ - m}\left( {{x_2}} \right){\text{e}^{ - {\rm{i}}{k_m}{x_1}}}.
\label{eq2.8}
\end{equation}

\noindent
and another incident wave of wavenumber $k_n$ propagating in the positive direction of axis $x_1$ in Eq.\eqref{eq2.2}.

Again, applying reciprocal theorem and considering the closed red line as the integral boundary shown in Fig., the following BIE can be obtained

\begin{equation}
\int\limits_{{\Gamma _{1 - 4}}} {\hat u_i^{ - m}\left( {{x_2}} \right){\text{e}^{ - {\rm{i}}{k_m}{x_1}}}\hat \sigma _{ij}^n\left( {{x_2}} \right){\text{e}^{{\rm{i}}{k_n}{x_1}}}{n_j}} \text{d}\Gamma \left( \bm{x} \right)\; - \int\limits_{{\Gamma _{1 - 4}}} {\hat u_i^n\left( {{x_2}} \right){\text{e}^{{\rm{i}}{k_n}{x_1}}}\hat \sigma _{ij}^{ - m}\left( {{x_2}} \right){\text{e}^{ - {\rm{i}}{k_m}{x_1}}}{n_j}} \text{d}\Gamma \left( \bm{x} \right)\; = 0
\label{eq2.9}
\end{equation}

Given the traction free boundary conditions of the upper and lower surfaces, the above equation can be simplified as

\begin{equation}
\left( {{\text{e}^{ - {\rm{i}}{k_m}{x_b}}}{\text{e}^{{\rm{i}}{k_n}{x_b}}} - {\text{e}^{ - {\rm{i}}{k_m}{x_a}}}{\text{e}^{{\rm{i}}{k_n}{x_a}}}} \right) \cdot {\int\limits_{-h_B}^{{h_A}} {\left[ {\hat u_3^{ - m}\left( {{x_2}} \right)\hat \sigma _{31}^n\left( {{x_2}} \right) - \hat u_3^n\left( {{x_2}} \right)\hat \sigma _{31}^{ - m}\left( {{x_2}} \right)} \right]\text{d}{x_2}} }  = 0
\label{eq2.10}
\end{equation}

It is also clearly found that when $m=n$, the product of the second term in Eq.\eqref{eq2.8} automatically satisfies and when $m\ne n$, we can get

\begin{equation}
{\int\limits_{-h_B}^{{h_A}} {\left[ {\hat u_3^{ - m}\left( {{x_2}} \right)\hat \sigma _{31}^n\left( {{x_2}} \right) - \hat u_3^n\left( {{x_2}} \right)\hat \sigma _{31}^{ - m}\left( {{x_2}} \right)} \right]\text{d}{x_2}} }  \equiv 0
\label{eq2.11}
\end{equation}

\noindent
since the variables $x_a$ and $x_b$ are both arbitrary.

Because 
\begin{equation}
\hat u_3^{-m}\left( {{x_2}} \right)=\hat u_3^{m}\left( {{x_2}} \right),\ \ \sigma _{31}^{-m}\left( {{x_2}} \right) = -\mu {\rm{i}}{k_{m}}\hat u_3^{m}\left( {{x_2}} \right)
\label{eq2.12}
\end{equation}

Eq.\eqref{eq2.12} can be also simplified as

\begin{equation}
\int\limits_{-h_B}^{{h_A}} { {\hat u_3^m\left( {{x_2}} \right)\hat u_3^n\left( {{x_2}} \right)\text{d}{x_2}} }  \equiv 0
\label{eq2.13}
\end{equation}

\subsection{Far-field conditions at the virtual boundaries}

The scattered displacement and traction in the far field are obtained by expanding the wave fields in a series of orthogonal modes. At a fixed frequency, the dispersion equation has a finite number of real roots corresponding to propagating modes, and an infinite number of complex roots corresponding to non-propagating modes. Since non-propagating waves geometrically attenuate in the propagating direction, the far-field displacement solution can be approximated by only the propagating waves. Therefore, when the truncated points are located far enough from the source regions, the far-field displacements and tractions at the virtual boundaries can be expressed as

\begin{equation}
\begin{aligned}
u_3^{\text{sca}+}\left( x_1,x_2\right) &= {\sum_{n=1}^N {\alpha_n^+ u_3^{n}\left( x_1,x_2\right)}},\ \ t_3^{\text{sca}+}\left( x_1,x_2\right) = {\sum_{n=1}^N {\alpha_n^+ t_3^{n}\left( x_1,x_2\right)}},\ \ x_1,x_2\in \Gamma_2, \\
u_3^{\text{sca}-}\left( x_1,x_2\right) &= {\sum_{n=1}^N {\alpha_n^- u_3^{-n}\left( x_1,x_2\right)}},\ \ t_3^{\text{sca}-}\left( x_1,x_2\right) = {\sum_{n=1}^N {\alpha_n^- t_3^{-n}\left( x_1,x_2\right)}},\ \ x_1,x_2\in \Gamma_1.
\end{aligned}
\label{eq2.14}
\end{equation}

The above equation can be further reduced to
\begin{equation}
\begin{aligned}
u_3^{\text{sca}+}\left( x_1,x_2\right) &= {\sum_{n=1}^N {\alpha_n^+ \hat u_3^n\left( {{x_2}} \right){\text{e}^{{\rm{i}}{k_n}{x_1}}} }},\ \ \sigma _{31}^{\text{sca}+}\left( x_1,x_2\right) = {\sum_{n=1}^N {\alpha_n^+ \hat \sigma _{31}^n\left( {{x_2}} \right) {\text{e}^{{\rm{i}}{k_n}{x_1}}} }},\ \ x_1,x_2\in \Gamma_2, \\
u_3^{\text{sca}-}\left( x_1,x_2\right) &= {\sum_{n=1}^N {\alpha_n^- \hat u_3^{n}\left( {{x_2}} \right){\text{e}^{-{\rm{i}}{k_n}{x_1}}} }},\ \ \sigma _{31}^{\text{sca}-}\left( x_1,x_2\right) = {\sum_{n=1}^N {\alpha_n^- \hat \sigma _{31}^{-n}\left( {{x_2}} \right) {\text{e}^{-{\rm{i}}{k_n}{x_1}}} }},\ \ x_1,x_2\in \Gamma_1.
\end{aligned}
\label{eq2.15}
\end{equation}

Multiplying both ends of the equals sign of the above equation by function $\hat u_3^m\left( x_2 \right)\ (m=1,2,...,N)$ and integrating them along $\Gamma _1$ or $\Gamma _2$, then we can achieve

\begin{equation}
\begin{gathered}
{\int_{\Gamma _1}u_3^{\text{sca}+}\left( x_1,x_2\right) \hat u_3^m\left( x_2 \right) \text{d}\Gamma _1} = {\sum_{n=1}^N \alpha_n^+ \cdot \int\limits_{-h_B}^{h_A}{\hat u_3^n\left( x_2 \right)\hat u_3^m\left( x_2 \right) \text{d} x_2} \cdot \text{e}^{\text{i} k_n a}}, \ \left( x_1=a\right) \\
{\int_{\Gamma _1}u_3^{\text{sca}-}\left( x_1,x_2\right) \hat u_3^m\left( x_2 \right) \text{d}\Gamma _2} = {\sum_{n=1}^N \alpha_n^- \cdot \int\limits_{-h_B}^{h_A}{\hat u_3^n\left( x_2 \right)\hat u_3^m\left( x_2 \right) \text{d} x_2} \cdot \text{e}^{\text{i} k_n a}}. \ \left( x_1=-a\right) \\
(m=1,2,...,N)
\end{gathered}
\label{eq2.16}
\end{equation}

Then, Eq.\ref{eq2.16} can be rewritten into the matrix form

\begin{equation}
\begin{gathered}
{\left[ {{G^ + }} \right]_{N \times 1}} = {\left[ R \right]_{N \times N}}{\left[ {{\alpha ^ + }} \right]_{N \times 1}} \\
{\left[ {{G^ - }} \right]_{N \times 1}} = {\left[ R \right]_{N \times N}}{\left[ {{\alpha ^ - }} \right]_{N \times 1}}
\end{gathered}
\label{eq2.17}
\end{equation}

\noindent
where,

\begin{equation}
\begin{gathered}
{\left[ {{G^ + }} \right]_{N \times 1}} = \left[ {\begin{array}{*{20}{c}}
{\int_{ - {h_B}}^{{h_A}} {u_3^{{\rm{sca}} + }\left( {a,{x_2}} \right)\hat u_3^1\left( {{x_2}} \right){\rm{d}}{x_2}} }\\
{\int_{ - {h_B}}^{{h_A}} {u_3^{{\rm{sca}} + }\left( {a,{x_2}} \right)\hat u_3^2\left( {{x_2}} \right){\rm{d}}{x_2}} }\\
 \vdots \\
{\int_{ - {h_B}}^{{h_A}} {u_3^{{\rm{sca}} + }\left( {a,{x_2}} \right)\hat u_3^N\left( {{x_2}} \right){\rm{d}}{x_2}} }
\end{array}} \right],\quad {\left[ {{G^ - }} \right]_{N \times 1}} = \left[ {\begin{array}{*{20}{c}}
{\int_{ - {h_B}}^{{h_A}} {u_3^{{\rm{sca}} - }\left( {a,{x_2}} \right)\hat u_3^1\left( {{x_2}} \right){\rm{d}}{x_2}} }\\
{\int_{ - {h_B}}^{{h_A}} {u_3^{{\rm{sca}} - }\left( {a,{x_2}} \right)\hat u_3^2\left( {{x_2}} \right){\rm{d}}{x_2}} }\\
 \vdots \\
{\int_{ - {h_B}}^{{h_A}} {u_3^{{\rm{sca}} - }\left( {a,{x_2}} \right)\hat u_3^N\left( {{x_2}} \right){\rm{d}}{x_2}} }
\end{array}} \right], \\
{\left[ {{\alpha ^ + }} \right]_{N \times 1}} = {\left[ {\begin{array}{*{20}{c}}
{\alpha _1^ + }&{\alpha _2^ + }& \cdots &{\alpha _N^ + }
\end{array}} \right]^T},\quad {\left[ {{\alpha ^ - }} \right]_{N \times 1}} = {\left[ {\begin{array}{*{20}{c}}
{\alpha _1^ - }&{\alpha _2^ - }& \cdots &{\alpha _N^ - }
\end{array}} \right]^T}, \\
{\left[ R \right]_{N \times N}} = \left[ {\begin{array}{*{20}{c}}
{{R_{11}}}&{{R_{12}}}& \cdots &{{R_{1N}}}\\
{{R_{21}}}&{{R_{22}}}& \cdots &{{R_{2N}}}\\
 \vdots & \vdots & \ddots & \vdots \\
{{R_{N1}}}&{{R_{N2}}}& \cdots &{{R_{NN}}}
\end{array}} \right],\quad {R_{mn}} = \int_{ - {h_B}}^{{h_A}} {\hat u_3^n\left( {{x_2}} \right)\hat u_3^m\left( {{x_2}} \right){\rm{d}}{x_2}}  \cdot {{\rm{e}}^{{\rm{i}}{k_n}a}}.
\end{gathered}
\label{eq2.18}
\end{equation}

Using the mode orthogonality derived from Eq.\ref{eq2.7} and Eq.\ref{eq2.13}, the matrix ${\left[ R \right]_{N \times N}}$ can be simplified as the diagonal matrix form

\begin{equation}
\begin{gathered}
{\left[ R \right]_{N \times N}} = \left[ {\begin{array}{*{20}{c}}
{{R_{11}}}&0& \cdots &0\\
0&{{R_{22}}}& \cdots &0\\
 \vdots & \vdots & \ddots & \vdots \\
0&0& \cdots &{{R_{NN}}}
\end{array}} \right]
\end{gathered}
\label{eq2.19}
\end{equation}

Then Eq.\ref{eq2.17} can be rewritten as

\begin{equation}
\begin{gathered}
{\left[ {{\alpha ^ + }} \right]_{N \times 1}} = \left[ R \right]_{N \times N}^{ - 1}{\left[ {{G^ + }} \right]_{N \times 1}} \\
{\left[ {{\alpha ^ - }} \right]_{N \times 1}} = \left[ R \right]_{N \times N}^{ - 1}{\left[ {{G^ - }} \right]_{N \times 1}}
\end{gathered}
\label{eq2.20}
\end{equation}

\noindent
where

\begin{equation}
\begin{gathered}
\left[ R \right]_{N \times N}^{ - 1} = \left[ {\begin{array}{*{20}{c}}
{\frac{1}{{{R_{11}}}}}&0& \cdots &0\\
0&{\frac{1}{{{R_{22}}}}}& \cdots &0\\
 \vdots & \vdots & \ddots & \vdots \\
0&0& \cdots &{\frac{1}{{{R_{NN}}}}}
\end{array}} \right]
\end{gathered}
\label{eq2.21}
\end{equation}

In the finite element analysis, the displacement and traction in element $e$ can be expressed as

\begin{equation}
\begin{gathered}
u_3^e = \sum\limits_{J = 1}^{{N_e}} {{N_J}u_3^J} ,\quad t_3^e = \sum\limits_{J = 1}^{{N_e}} {{N_J}t_3^J} .
\end{gathered}
\label{eq2.22}
\end{equation}

\noindent
where $N_J$ is the shape function at node $J$ and $u_3^J$ is the displacement at node $J$ in element $e$.

Therefore, the matrix ${\left[ {{G^ + }} \right]_{N \times 1}}$ and ${\left[ {{G^ - }} \right]_{N \times 1}}$ can be derived as

\begin{equation}
\begin{gathered}
{\left[ {{G^ + }} \right]_{N \times 1}} = \left[ {\begin{array}{*{20}{c}}
{\sum\limits_{e \in {\Gamma _1}} {\sum\limits_{J = 1}^{{N_e}} {\int\limits_{{\Gamma _e}} {{N_J}\hat u_3^1\left( {{x_2}} \right){\rm{d}}{\Gamma _e}} u_3^J} } }\\
{\sum\limits_{e \in {\Gamma _1}} {\sum\limits_{J = 1}^{{N_e}} {\int\limits_{{\Gamma _e}} {{N_J}\hat u_3^1\left( {{x_2}} \right){\rm{d}}{\Gamma _e}} u_3^J} } }\\
 \vdots \\
{\sum\limits_{e \in {\Gamma _1}} {\sum\limits_{J = 1}^{{N_e}} {\int\limits_{{\Gamma _e}} {{N_J}\hat u_3^1\left( {{x_2}} \right){\rm{d}}{\Gamma _e}} u_3^J} } }
\end{array}} \right] = {\left[ {{{\hat G}^ + }} \right]_{N \times P}}{\left[ U \right]_{P \times 1}}
\end{gathered}
\label{eq2.23}
\end{equation}

\noindent
and

\begin{equation}
\begin{gathered}
{\left[ {{G^ - }} \right]_{N \times 1}} = \left[ {\begin{array}{*{20}{c}}
{\sum\limits_{e \in {\Gamma _2}} {\sum\limits_{J = 1}^{{N_e}} {\int\limits_{{\Gamma _e}} {{N_J}\hat u_3^1\left( {{x_2}} \right){\rm{d}}{\Gamma _e}} u_3^J} } }\\
{\sum\limits_{e \in {\Gamma _2}} {\sum\limits_{J = 1}^{{N_e}} {\int\limits_{{\Gamma _e}} {{N_J}\hat u_3^1\left( {{x_2}} \right){\rm{d}}{\Gamma _e}} u_3^J} } }\\
 \vdots \\
{\sum\limits_{e \in {\Gamma _2}} {\sum\limits_{J = 1}^{{N_e}} {\int\limits_{{\Gamma _e}} {{N_J}\hat u_3^1\left( {{x_2}} \right){\rm{d}}{\Gamma _e}} u_3^J} } }
\end{array}} \right] = {\left[ {{{\hat G}^ - }} \right]_{N \times P}}{\left[ U \right]_{P \times 1}}
\end{gathered}
\label{eq2.24}
\end{equation}

\noindent
where matrices ${\left[ {{{\hat G}^ + }} \right]_{N \times P}}$ and ${\left[ {{{\hat G}^ + }} \right]_{N \times P}}$ are assembled from known element ${\int\limits_{{\Gamma _e}} {{N_J}\hat u_3^1\left( {{x_2}} \right){\rm{d}}{\Gamma _e}} }$, $P$ is the total node number and the matrix ${\left[ U \right]_{P \times 1}}$ is the node displacement matrix in global coordinates.

Due to Eq.\ref{eq2.23} and Eq.\ref{eq2.24}, the unknown scattered coefficient matrices in Eq.\ref{eq2.20} can be rewritten as

\begin{equation}
\begin{gathered}
\begin{array}{l}
{\left[ {{\alpha ^ + }} \right]_{N \times 1}} = \left[ R \right]_{N \times N}^{ - 1}{\left[ {{{\hat G}^ + }} \right]_{N \times P}}{\left[ U \right]_{P \times 1}},\\
{\left[ {{\alpha ^ - }} \right]_{N \times 1}} = \left[ R \right]_{N \times N}^{ - 1}{\left[ {{{\hat G}^ - }} \right]_{N \times P}}{\left[ U \right]_{P \times 1}}.
\end{array}
\end{gathered}
\label{eq2.25}
\end{equation}

According to Eq.\ref{eq2.14} and Eq.\ref{eq2.15}, tractions at the virtual boundaries can be expressed as the matrix forms

\begin{equation}
\begin{gathered}
\begin{array}{l}
t_3^{{\rm{sca}} + }\left( {{x_1},{x_2}} \right) = {\left[ {{\sigma ^ + }} \right]_{1 \times N}}{\left[ {{\alpha ^ + }} \right]_{N \times 1}},\\
t_3^{{\rm{sca}} - }\left( {{x_1},{x_2}} \right) = {\left[ {{\sigma ^ - }} \right]_{1 \times N}}{\left[ {{\alpha ^ - }} \right]_{N \times 1}}.
\end{array}
\end{gathered}
\label{eq2.26}
\end{equation}

\noindent
where

\begin{equation}
\begin{gathered}
{\left[ {{\sigma ^ + }} \right]_{1 \times N}} = {\left[ {\begin{array}{*{20}{c}}
{\mu {\rm{i}}{k_1}\hat u_3^1\left( {{x_2}} \right){{\rm{e}}^{{\rm{i}}{k_1}a}}}\\
{\mu {\rm{i}}{k_2}\hat u_3^2\left( {{x_2}} \right){{\rm{e}}^{{\rm{i}}{k_2}a}}}\\
 \vdots \\
{\mu {\rm{i}}{k_N}\hat u_3^N\left( {{x_2}} \right){{\rm{e}}^{{\rm{i}}{k_N}a}}}
\end{array}} \right]^T}
\end{gathered}
\label{eq2.27}
\end{equation}

\begin{equation}
\begin{gathered}
{\left[ {{\sigma ^ - }} \right]_{1 \times N}} = {\left[ {\begin{array}{*{20}{c}}
{-\mu {\rm{i}}{k_1}\hat u_3^1\left( {{x_2}} \right){{\rm{e}}^{{\rm{i}}{k_1}a}}}\\
{-\mu {\rm{i}}{k_2}\hat u_3^2\left( {{x_2}} \right){{\rm{e}}^{{\rm{i}}{k_2}a}}}\\
 \vdots \\
{-\mu {\rm{i}}{k_N}\hat u_3^N\left( {{x_2}} \right){{\rm{e}}^{{\rm{i}}{k_N}a}}}
\end{array}} \right]^T}
\end{gathered}
\label{eq2.28}
\end{equation}

Using conditions in Eq.\ref{eq2.25}, Eq.\ref{eq2.26} can be rewritten as

\begin{equation}
\begin{gathered}
\begin{array}{l}
t_3^{{\rm{sca}} + }\left( {{x_1},{x_2}} \right) = {\left[ {{\sigma ^ + }} \right]_{1 \times N}}\left[ R \right]_{N \times N}^{ - 1}{\left[ {{{\hat G}^ + }} \right]_{N \times P}}{\left[ U \right]_{P \times 1}},\\
t_3^{{\rm{sca}} - }\left( {{x_1},{x_2}} \right) = {\left[ {{\sigma ^ - }} \right]_{1 \times N}}\left[ R \right]_{N \times N}^{ - 1}{\left[ {{{\hat G}^ - }} \right]_{N \times P}}{\left[ U \right]_{P \times 1}}.
\end{array}
\end{gathered}
\label{eq2.29}
\end{equation}

The equivalent nodal force at node $J$ of element $e$ can be expressed as

\begin{equation}
\begin{array}{l}
\int\limits_{{\Gamma _e}} {t_3^{{\rm{sca}} + }\left( {{x_1},{x_2}} \right){N_J}{\rm{d}}{\Gamma _e}}  = {\left[ {F_e^ + } \right]_{1 \times N}}\left[ R \right]_{N \times N}^{ - 1}{\left[ {{{\hat G}^ + }} \right]_{N \times P}}{\left[ U \right]_{P \times 1}},\\
\int\limits_{{\Gamma _e}} {t_3^{{\rm{sca}} - }\left( {{x_1},{x_2}} \right){N_J}{\rm{d}}{\Gamma _e}}  = {\left[ {F_e^ - } \right]_{1 \times N}}\left[ R \right]_{N \times N}^{ - 1}{\left[ {{{\hat G}^ - }} \right]_{N \times P}}{\left[ U \right]_{P \times 1}}.
\end{array}
\label{eq2.30}
\end{equation}

\noindent
where

\begin{equation}
\begin{gathered}
{\left[ {{F_e^+ }} \right]_{1 \times N}} = {\left[ {\begin{array}{*{20}{c}}
{\int\limits_{{\Gamma _e}} {\mu {\rm{i}}{k_1}\hat u_3^1\left( {{x_2}} \right){{\rm{e}}^{{\rm{i}}{k_1}a}}{N_J}{\rm{d}}{\Gamma _e}} }\\
{\int\limits_{{\Gamma _e}} {\mu {\rm{i}}{k_2}\hat u_3^2\left( {{x_2}} \right){{\rm{e}}^{{\rm{i}}{k_2}a}}{N_J}{\rm{d}}{\Gamma _e}} }\\
 \vdots \\
{\int\limits_{{\Gamma _e}} {\mu {\rm{i}}{k_N}\hat u_3^N\left( {{x_2}} \right){{\rm{e}}^{{\rm{i}}{k_N}a}}{N_J}{\rm{d}}{\Gamma _e}} }
\end{array}} \right]^T}
\end{gathered}
\label{eq2.31}
\end{equation}

\begin{equation}
\begin{gathered}
{\left[ {{F_e^- }} \right]_{1 \times N}} = {\left[ {\begin{array}{*{20}{c}}
{-\int\limits_{{\Gamma _e}} {\mu {\rm{i}}{k_1}\hat u_3^1\left( {{x_2}} \right){{\rm{e}}^{{\rm{i}}{k_1}a}}{N_J}{\rm{d}}{\Gamma _e}} }\\
{-\int\limits_{{\Gamma _e}} {\mu {\rm{i}}{k_2}\hat u_3^2\left( {{x_2}} \right){{\rm{e}}^{{\rm{i}}{k_2}a}}{N_J}{\rm{d}}{\Gamma _e}} }\\
 \vdots \\
{-\int\limits_{{\Gamma _e}} {\mu {\rm{i}}{k_N}\hat u_3^N\left( {{x_2}} \right){{\rm{e}}^{{\rm{i}}{k_N}a}}{N_J}{\rm{d}}{\Gamma _e}} }
\end{array}} \right]^T}
\end{gathered}
\label{eq2.32}
\end{equation}

It should be noted that matrices in Eq.\ref{eq2.31} and Eq.\ref{eq2.32} are known. Then, we assembly local matrices ${\left[ {F_e^ + } \right]_{1 \times N}}$ and ${\left[ {F_e^ - } \right]_{1 \times N}}$ into global matrices ${\left[ {F_g^ + } \right]_{P \times N}}$ and ${\left[ {F_g^ - } \right]_{P \times N}}$ respectively. 
Therefore, the global equivalent nodal force matrix can be expressed as
\begin{equation}
\begin{aligned}
{\left[ {{F^{{\rm{sca}}}}} \right]_{P \times 1}} &= {\left[ {F_g^ + } \right]_{P \times N}}\left[ R \right]_{N \times N}^{ - 1}{\left[ {{{\hat G}^ + }} \right]_{N \times P}}{\left[ U \right]_{P \times 1}} + {\left[ {F_g^ - } \right]_{P \times N}}\left[ R \right]_{N \times N}^{ - 1}{\left[ {{{\hat G}^ - }} \right]_{N \times P}}{\left[ U \right]_{P \times 1}} \\
&= \left\{ {{{\left[ {F_g^ + } \right]}_{P \times N}}\left[ R \right]_{N \times N}^{ - 1}{{\left[ {{{\hat G}^ + }} \right]}_{N \times P}} + {{\left[ {F_g^ - } \right]}_{P \times N}}\left[ R \right]_{N \times N}^{ - 1}{{\left[ {{{\hat G}^ - }} \right]}_{N \times P}}} \right\}{\left[ U \right]_{P \times 1}} \\
&= {\left[ {{{\bar F}^{{\rm{sca}}}}} \right]_{P \times P}}{\left[ U \right]_{P \times 1}}
\end{aligned}
\label{eq2.33}
\end{equation}

\subsection{DtN-FEM formulation}

The elasto-dynamic finite element formulation is 
\begin{equation}
\begin{gathered}
\left( {\left[ K \right] - {\omega ^2}\left[ M \right]} \right)\left[ {{U^{{\rm{tot}}}}} \right] = \left[ {{F^{{\rm{tot}}}}} \right]
\end{gathered}
\label{eq2.34}
\end{equation}

\noindent
where $\omega $ is the angular frequency, 

\begin{equation}
\begin{gathered}
\left[ K \right] = \int\limits_S {{{\left[ B \right]}^T}\left[ D \right]\left[ B \right]{\rm{d}}S} ,\quad \left[ M \right] = \rho \int\limits_S {{{\left[ N \right]}^T}\left[ N \right]{\rm{d}}S} , \\
\left[ {{U^{{\rm{tot}}}}} \right] = \left[ {{U^{{\rm{inc}}}}} \right] + \left[ {{U^{{\rm{sca}}}}} \right],\quad \left[ {{F^{{\rm{tot}}}}} \right] = \left[ {{F^{{\rm{inc}}}}} \right] + \left[ {{F^{{\rm{sca}}}}} \right].
\end{gathered}
\label{eq2.35}
\end{equation}

\noindent
and ${\left[ N \right]}$ is the matrix of shape function, ${\left[ B \right]}$ is the geometric matrix, ${\left[ D \right]}$ is the elastic matrix, ${\left[ K \right]}$ is the stiffness matrix, ${\left[ M \right]}$ is the mass matrix, and $\rho $ is the material density.

Utilizing Eq.\ref{eq2.33}, Eq.\ref{eq2.34} can be finally simplified as
\begin{equation}
\begin{gathered}
\left[ {{U^{{\rm{sca}}}}} \right] = {\left\{ {\left( {\left[ K \right] - {\omega ^2}\left[ M \right]} \right) - \left[ {{{\bar F}^{{\rm{sca}}}}} \right]} \right\}^{ - 1}}\left\{ {\left[ {{F^{{\rm{inc}}}}} \right] - \left( {\left[ K \right] - {\omega ^2}\left[ M \right]} \right)\left[ {{U^{{\rm{inc}}}}} \right]} \right\}
\end{gathered}
\label{eq2.35}
\end{equation}

\noindent
which is the basic DtN-FEM formulation.

\section{Numerical examples}

In this section, the validity will be presented in order to verify the correctness and efficiency of our proposed DtN-FEM. Here, the DtN-FEM results will be compared with the BEM results, and the error analysis of energy will also be demonstrated. Furthermore, parametric analysis about the effect on different materials, various lengths of the delamination and location of the interface will be discussed detailedly.

\begin{figure}[h]
	\centering
	\includegraphics[scale=0.6]{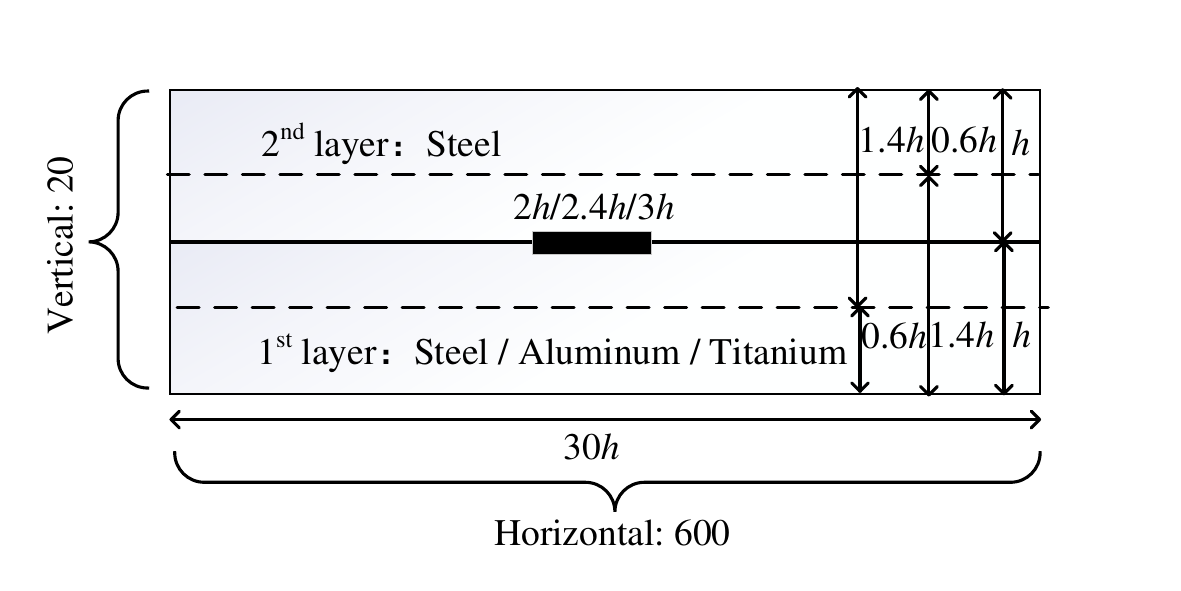}
	\caption{Model diagram for parametric studies}
	\label{fig8}
 \end{figure}

In the following numerical examples, an eight-node quadrilateral isoparametric element will be adopted. In Fig.\ref{fig8}, the element number for horizontal direction and vertical direction are shown, and other parameters for numerical calculation are also demonstrated. To satisfy the far-field assumption of SH guided wave, the horizontal length of the plate should be long enough. The total thickness of the plate is 1mm. The length of the plate $L$ is equal to $30h$, and the thickness of the plate is equal to $2h$, where $h$ is the half-width of the plate. It is worthy to mention that the relative ratio of the thickness of the upper and lower layers will change, but the total thickness of the plate is equal to 1 mm. The material properties of each plate are given in Tab.\ref{tab1}.

\begin{table}[h]
    \caption{Material constants}    
    \centering
    \begin{tabular}{cccc}
            \hline
            Material & $\lambda$(GPa) & $\mu$(GPa) & $\rho$(g/c$\rm{m}^3$)  \\  
            \hline
            Steel & 115.5 & 79.0 & 7.8 \\ 
            Aluminum & 58.2 & 26.1 & 2.7 \\
            Titanium & 66.9 & 44.6 & 4.5 \\
            \hline
            \end{tabular}

    \label{tab1}
\end{table}

\subsection{Validity of of our proposed DtN FEM}

In order to verify the correctness and accuracy of our proposed DtN FEM, normalized displacements on both top and bottom boundaries are plotted. Here, the normalized displacement is defined as

\begin{equation}
\frac{u_{3}^{- sca}}{\sum_{n = 1}^{N}{\alpha^{-}_nu_{3}^{- n}}},~~~~~~\frac{u_{3}^{+ sca}}{\sum_{n = 1}^{N}{\alpha^{+}_nu_{3}^{+ n}}}.
\label{eq3.18}
\end{equation}

\begin{figure}[h]
	\centering
	\includegraphics[scale=0.55]{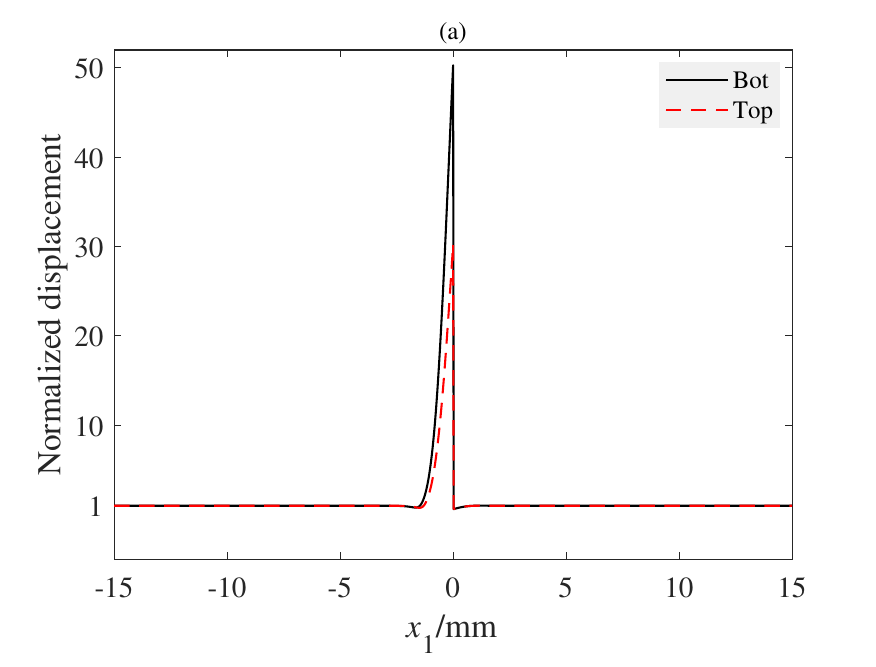}
    \includegraphics[scale=0.55]{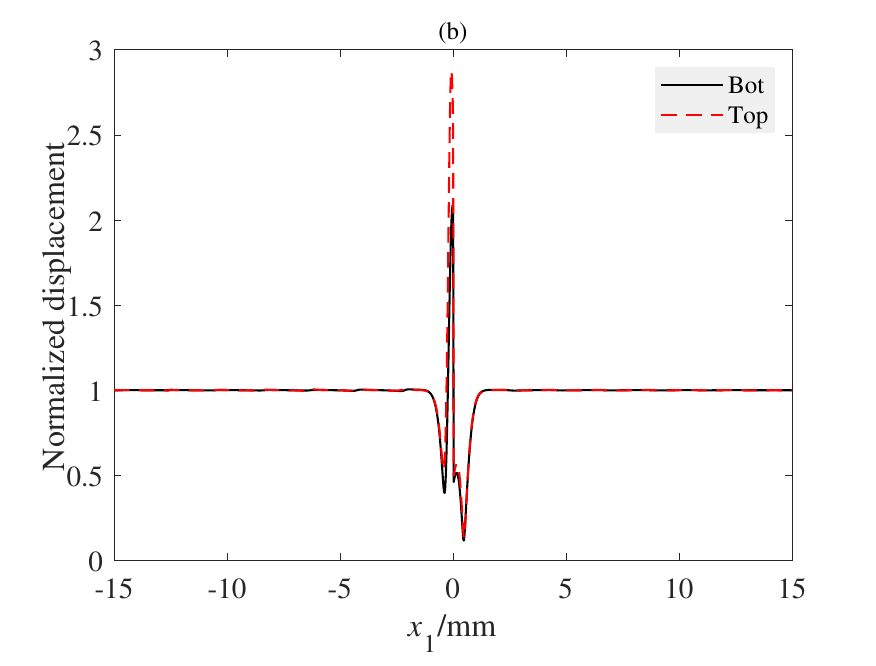}
	\caption{The normalized scattered displacement on the top and bottom boundaries in a 1mm thick aluminum-steel plate: (a) frequency is 2MHz; (b) frequency is 5MHz.}
	\label{fig9}
 \end{figure}

In this case, two different excitation frequencies 2MHz and 5MHz are selected, and the lowest guided mode is chosen as the incident wave. The normalized scattered displacements on the top and bottom boundaries in a 1mm thick aluminum-steel plate are shown in Fig.\ref{fig9}. It is clearly found that normalized displacements at both ends converges to 1, which indicates that as the propagation distance increases, the non-propagating modes decay rapidly, leaving only the guided wave modes. Therefore, this scattering phenomenon coincides perfectly with the previously defined far-field assumption which verifies the correctness and accuracy of our proposed DtN FEM.

on the other hand, we use the formula of energy conservation to verify the correctness and accuracy. The error analysis criterion is

\begin{equation}
\varepsilon =  {1 - {\sum\limits_{n = 1}^{N_{p}}\frac{E_{n}^{Re} + E_{n}^{Tr}}{E_{n}^{Inv}}}} 
\label{eq3.19}
\end{equation}

\noindent
where $E_n^{Re}$ indicates the reflection energy, $E_n^{Tr}$ indicates the transmission energy, $E_n^{Inv}$ is the energy of incident wave.

\begin{figure}[p]
\centering
 \subfigure[Different materials where the material of the lower layer is fixed at steel, while the material of the upper layer can be steel/aluminum/titanium, St, Al and Ti are abbreviations respectively of steel, aluminum and titanium.]{\includegraphics[scale=0.55]{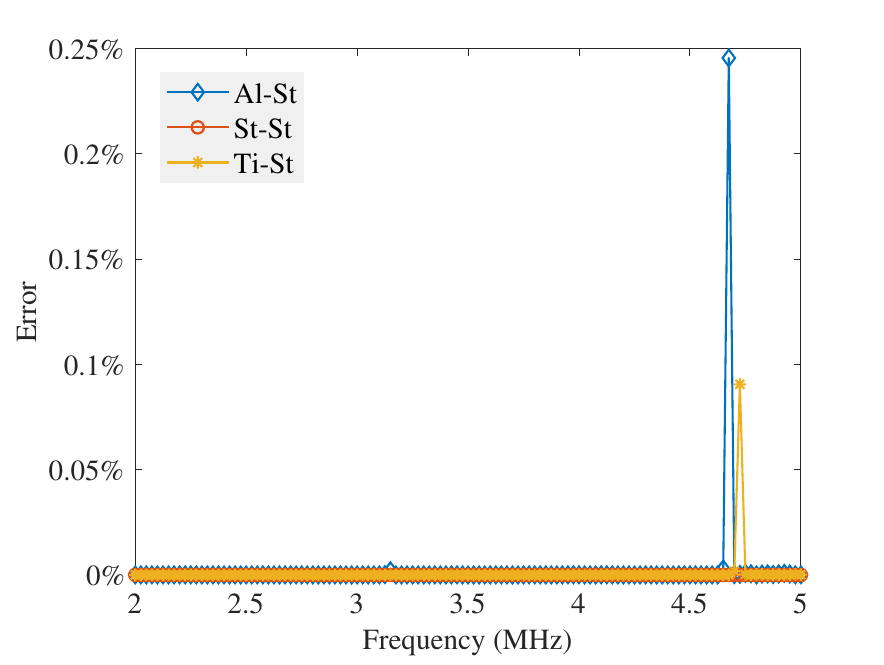}}
 \quad
\subfigure[Various lengths of the delamination where various lengths of delamination ($h$, $1.2h$ and $1.5h$) are considered, an aluminum-steel plate is adopted as the model, the thickness of each layer is assumed as $0.5h$, and the second guided mode is selected as the incident wave.]{\includegraphics[scale=0.55]{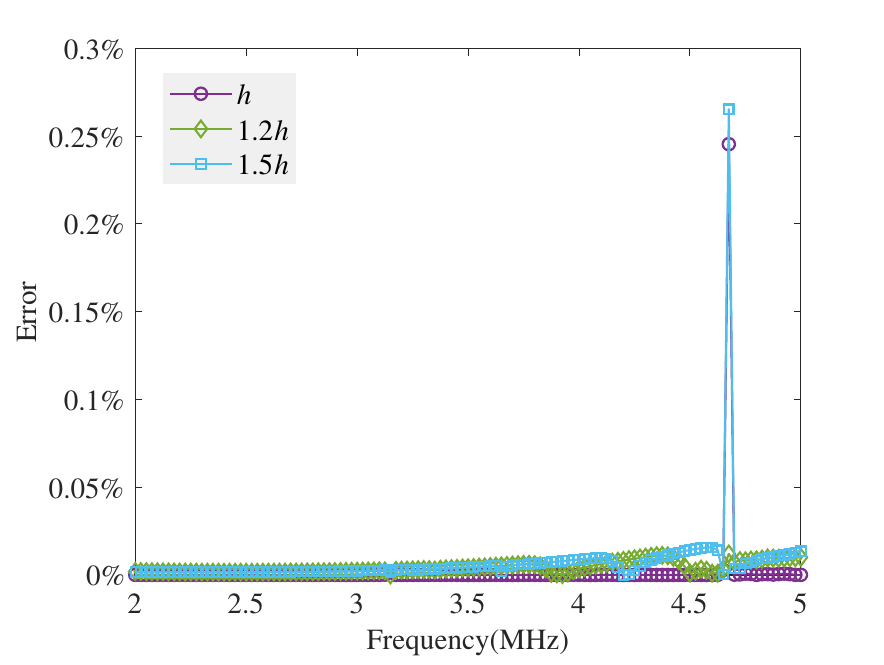}}
 \quad
\subfigure[Different location of the interface where different ratios of the first layer thickness and the total thickness (0.3, 0.5 and 0.7) are considered, an aluminum-steel plate is adopted as the model, the total thickness is assumed as $h$, and the second guided mode is selected as the incident wave.]{\includegraphics[scale=0.55]{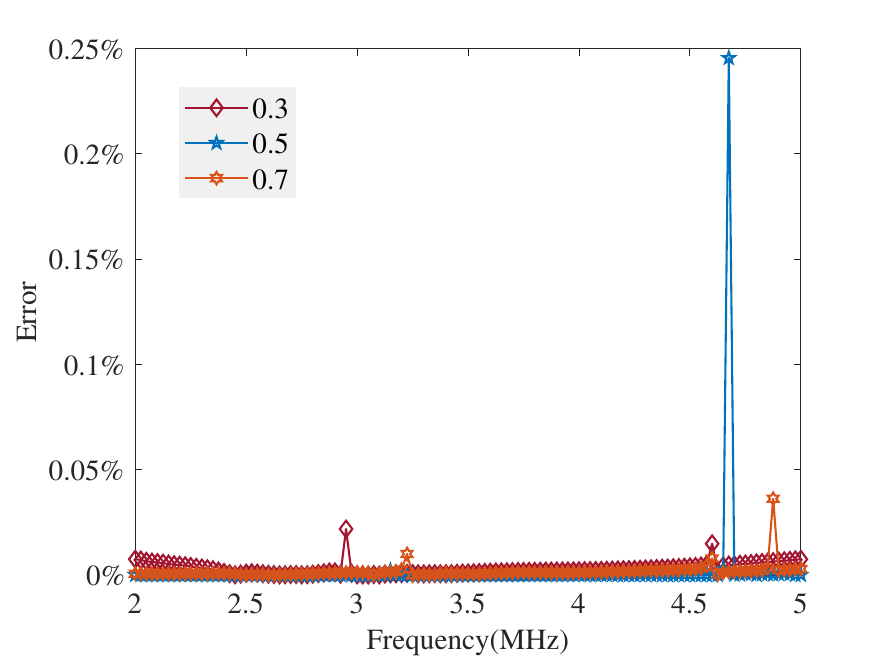}}

\caption{Error analysis of energy balance}
\label{fig10}
\end{figure}

It is clear from Fig.\ref{fig10} that the maximum error of energy balance is only around 0.25$\%$ and the average error is lower than $0.05\%$ which ensures sufficient calculation accuracy. Furthermore, the extreme value of the error appears when the material of two layers or the delamination length is changed, which is more obvious in high frequency domain.
 
\begin{figure}[h]
	\centering
	\includegraphics[scale=0.55]{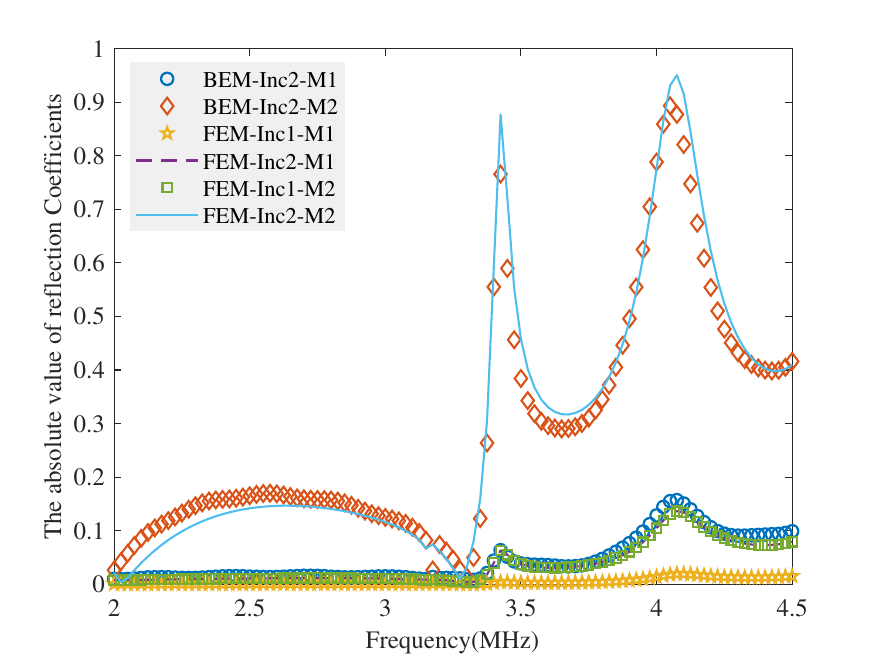}
    \includegraphics[scale=0.55]{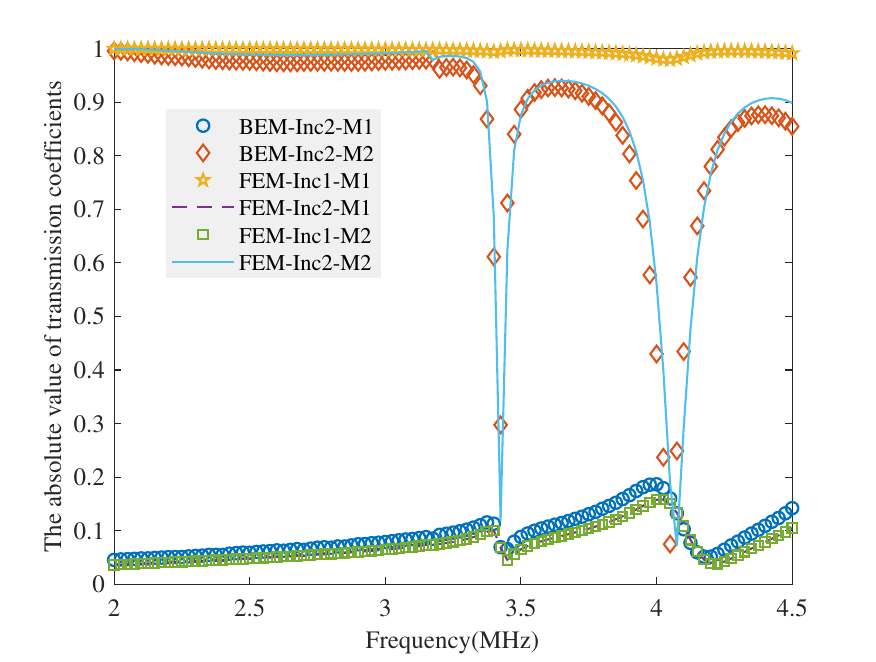}
	\caption{The absolute value of reflection and transmission coefficients calculted by proposed DtN-FEM and BEM}
	\label{fig12}
 \end{figure}

\begin{figure}[h]
	\centering
	\includegraphics[scale=0.7]{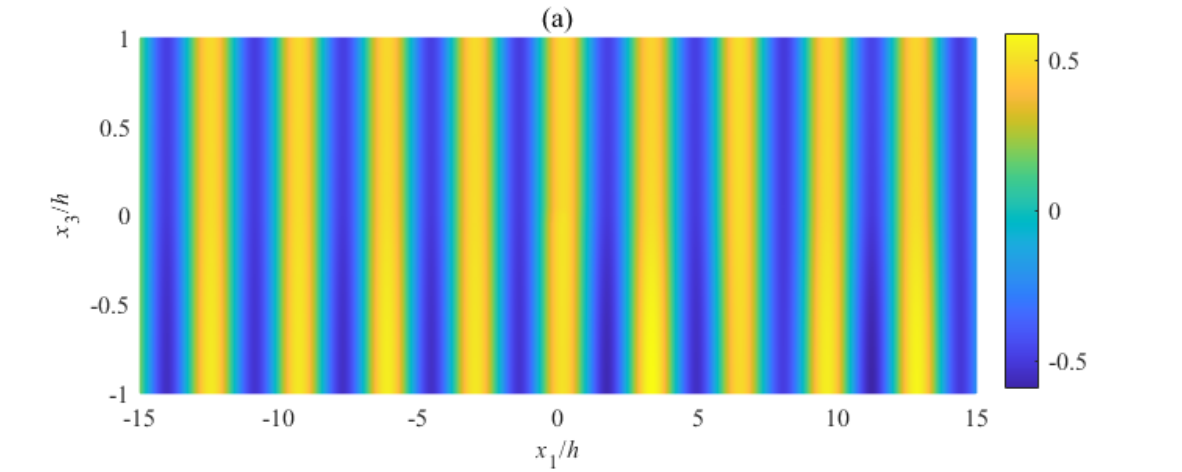}
    \includegraphics[scale=0.7]{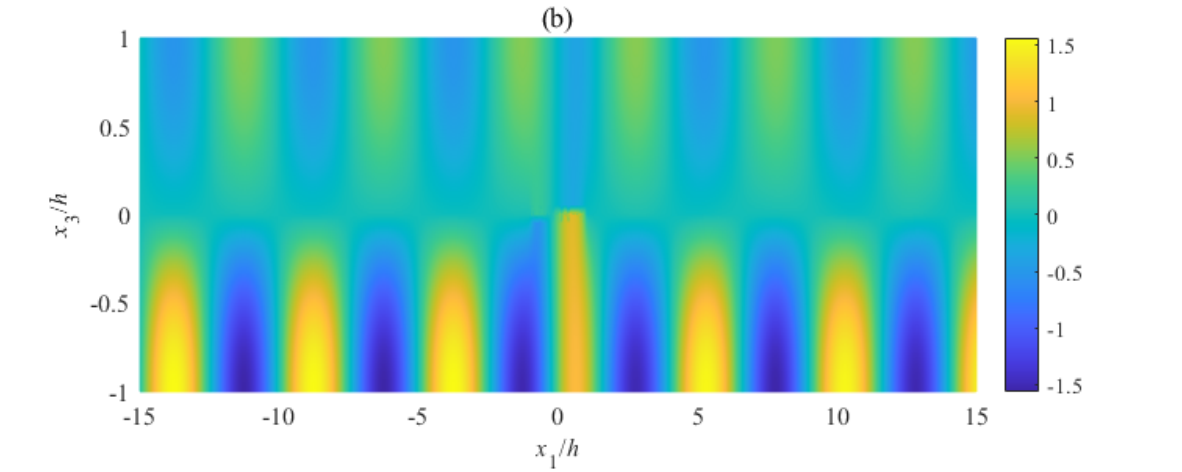}
	\caption{Total field displacement on a 1 mm thick aluminum-steel plate calculated by proposed DtN-FEM FEM when frequency is 2MHz: (a) The first guided mode is selected as the incident wave; (b) The second guided mode is selected as the incident wave.}
	\label{fig13}
 \end{figure}

In the end of this section, we compare the numerical results of our proposed DtN FEM with BEM's results from Wongthongsiri and Hirose [22]. We can definitely observe from Fig.\ref{fig12} that 
the reflection and transmission coefficients computed by FEM are perfectly consistent with BEM's results which verifies the correctness and accuracy of our proposed DtN FEM. Besides, it should be noted that when the first guided mode is selected as the incident wave, the absolute value of reflection coefficient on mode 1 is really small compared to that of transmission coefficient, which indicates that the first guided wave almost transmits through the debonding region. Total field displacements on a 1 mm thick aluminum-steel plate calculated by proposed DtN FEM are plotted in Fig.\ref{fig13} where different guided mode is selected as the incident wave. In Fig.\ref{fig13}(a), almost no scattering occurs since the first guided mode accounts for the main component. While the obvious scattering can be found in Fig.\ref{fig13}(b) since the second guided mode accounts for the main component.

\subsection{Parametric analysis}

\subsubsection{Different materials}

\begin{figure}[h]
	\centering
	\includegraphics[scale=0.55]{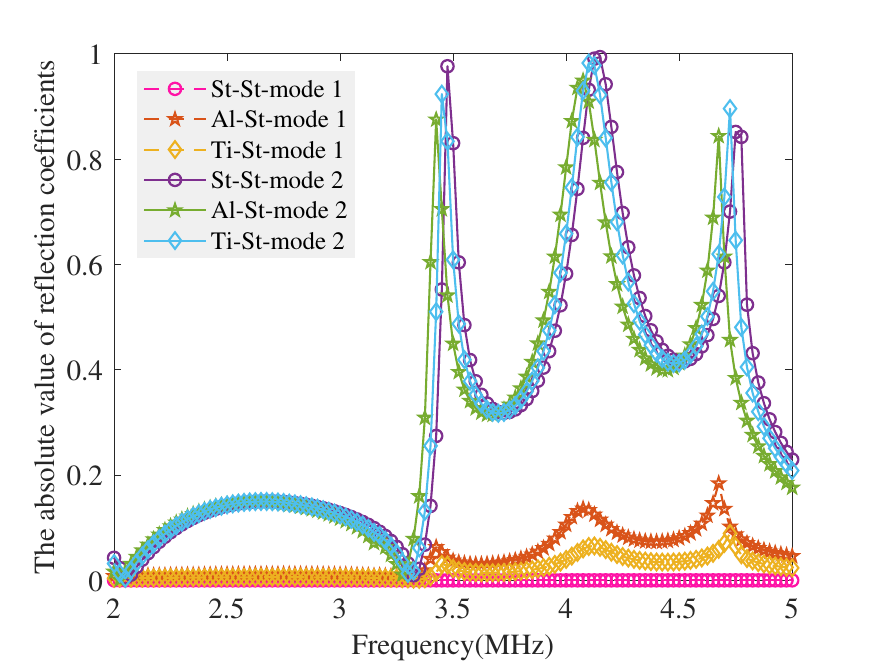}
    \includegraphics[scale=0.55]{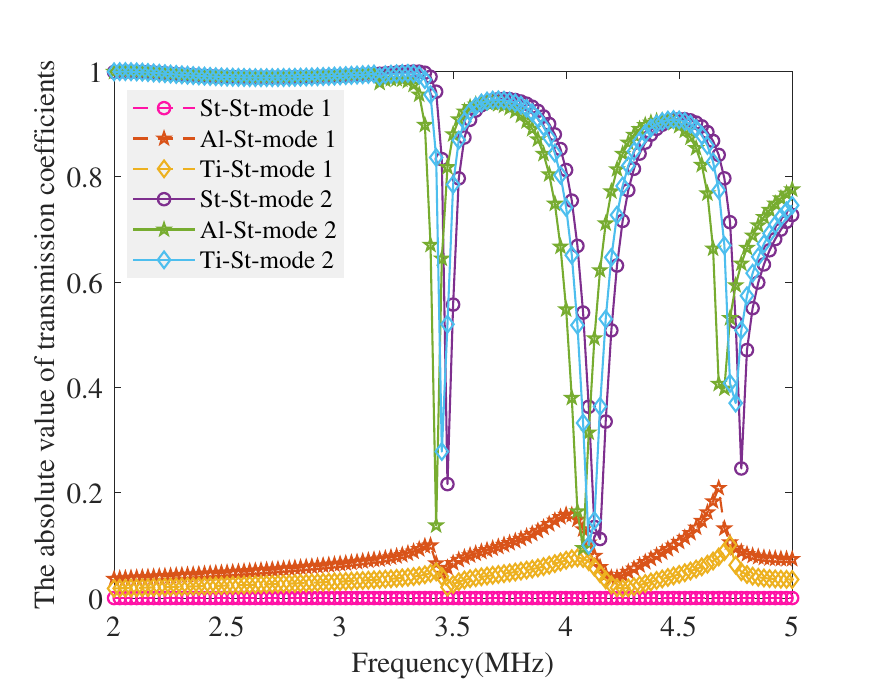}
	\caption{The absolute value of reflection and transmission coefficients varies with the frequency where St, Al and Ti are abbreviations respectively of steel, aluminum and titanium, the thickness of each layer is assumed as $0.5h$, and the second guided mode is selected as the incident wave}
	\label{fig14}
 \end{figure}

\begin{figure}[h]
	\centering
	\includegraphics[scale=0.7]{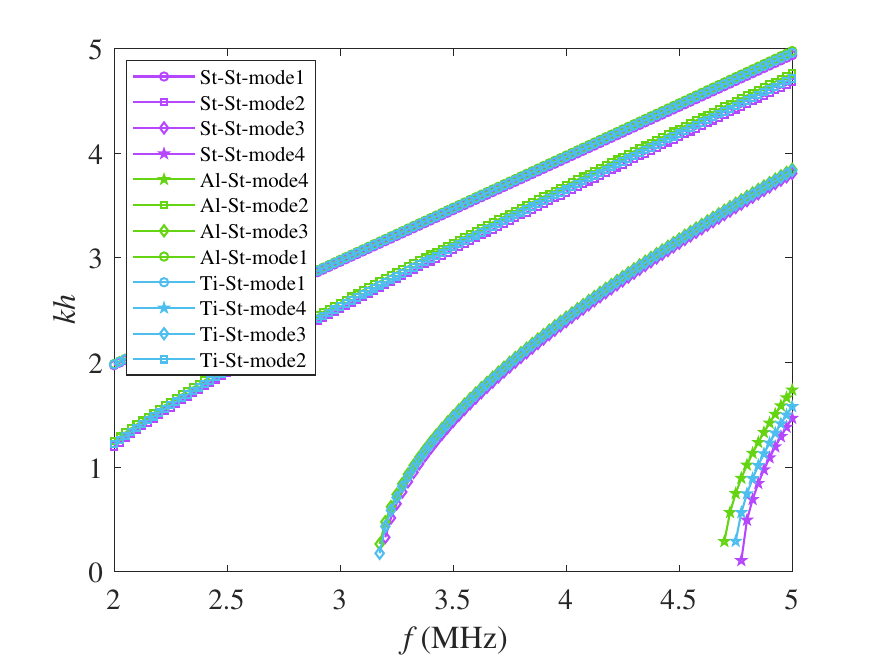}
	\caption{Dispersion curves due to different materials}
	\label{figmaterial_disper}
 \end{figure}

\begin{figure}[h]
\centering
 \subfigure[St-St]{\includegraphics[scale=0.65]{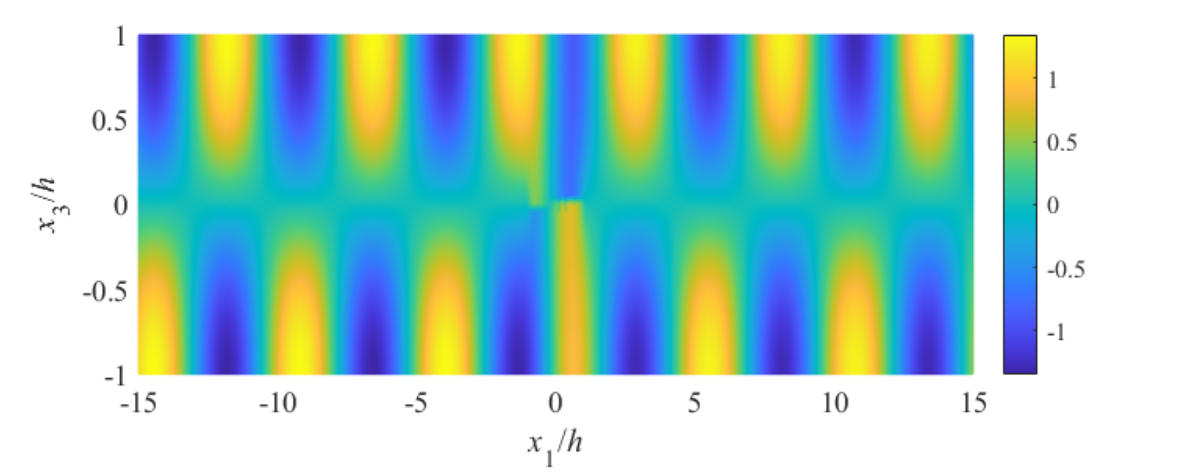}}
 \quad
\subfigure[Ti-St]{\includegraphics[scale=0.65]{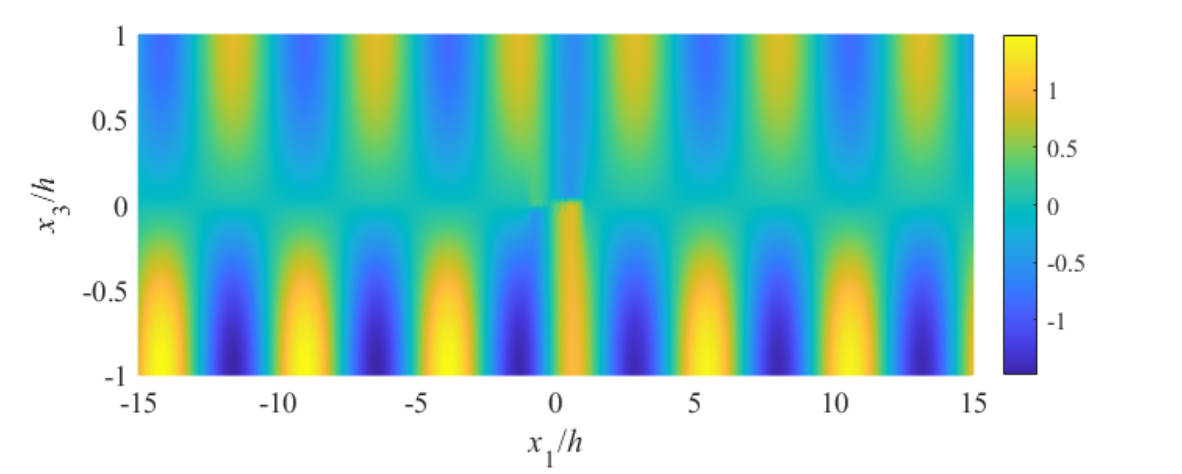}}
\caption{The total field displacement where the frequency is fixed at 2 MHz, the thickness of each layer is assumed as $0.5h$, and the second guided mode is selected as the incident wave.}
\label{fig15}
\end{figure}

In this section, the material of the lower layer is fixed at steel, while the material of the upper layer can be steel/aluminum/titanium. 
The thickness of each layer is assumed as $0.5h$, and the second guided mode is selected as the incident wave. Fig.\ref{fig14} demonstrates that the absolute value of reflection and transmission coefficients varies with the frequency where St, Al and Ti are abbreviations respectively of steel, aluminum and titanium. It is clearly observed from Fig.\ref{fig14} that all curves have a large fluctuation when the new guided modes appear due to dispersion curves Fig.\ref{figmaterial_disper}, and the absolute value of reflection and transmission coefficients due to mode 2 is larger compared to mode 1 which indicates that the second guided mode accounts for the main component. It should also be noted that the new mode of case Al-St is the first to appear and 
both the absolute values of reflection and transmission coefficients of mode 1 due to case Al-St are the largest compared with other two cases. It can be observed clearly from the Fig.\ref{fig15} that the total displacement of case St-St is almost anti-symmetric about the interface since the anti-symmetric mode 2 accounts for the main component, and the geometry and material property are both symmetric about the interface.

\subsubsection{Various lengths of the delamination}

\begin{figure}[h]
	\centering
	\includegraphics[scale=0.55]{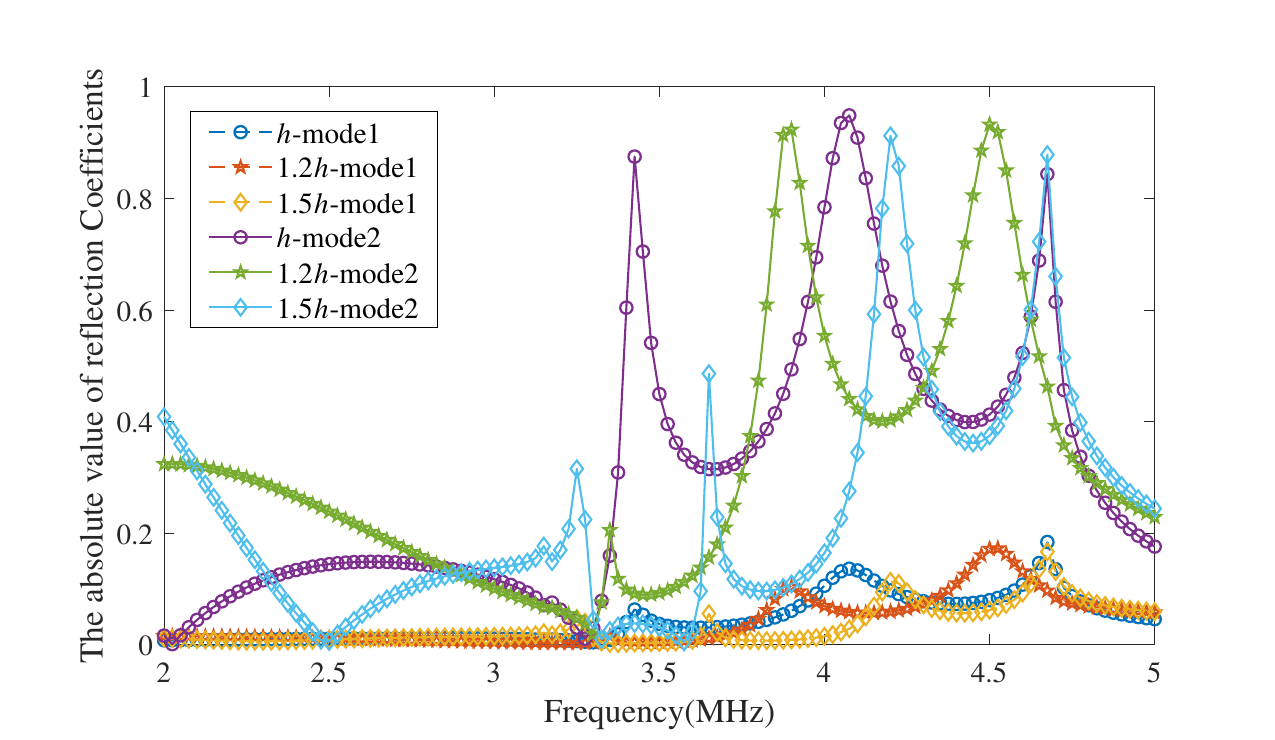}
    \includegraphics[scale=0.55]{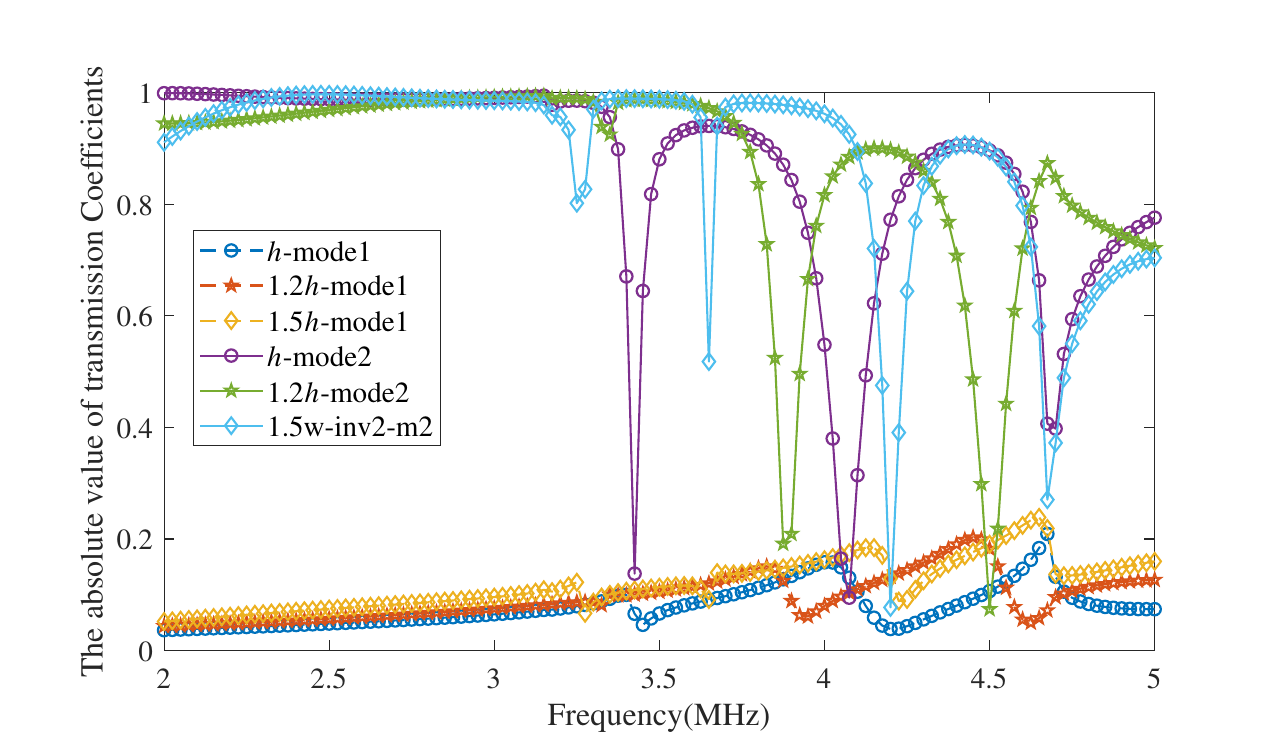}
	\caption{The absolute value of reflection and transmission coefficients varies with different length of delamination ($h$, $1.2h$ and $1.5h$) where St, Al and Ti are abbreviations respectively of steel, aluminum and titanium, the thickness of each layer is assumed as $0.5h$, and the second guided mode is selected as the incident wave}
	\label{fig16}
 \end{figure}

\begin{figure}[h]
\centering
 \subfigure[Length of delamination: $1.2h$]{\includegraphics[scale=0.7]{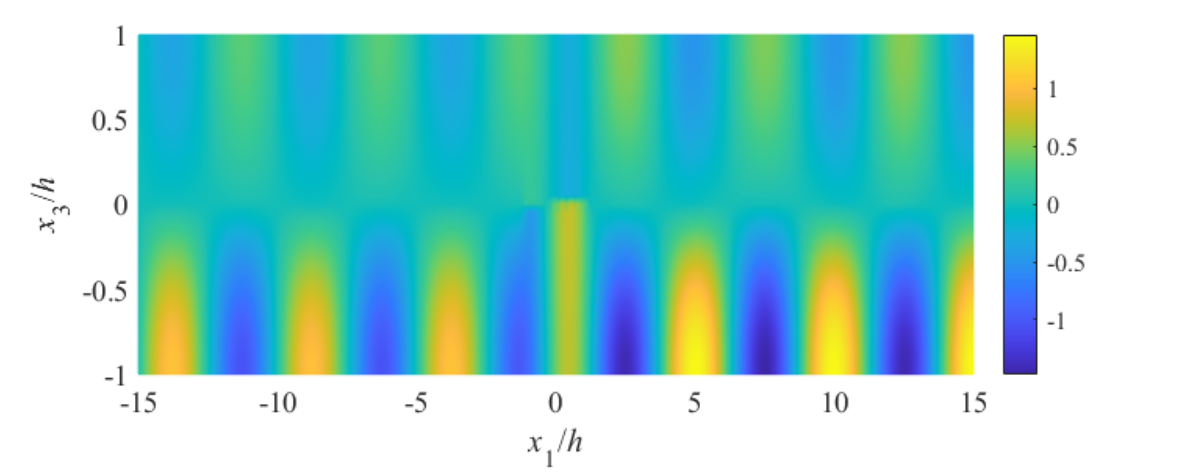}}
 \quad
\subfigure[Length of delamination: $1.5h$]{\includegraphics[scale=0.7]{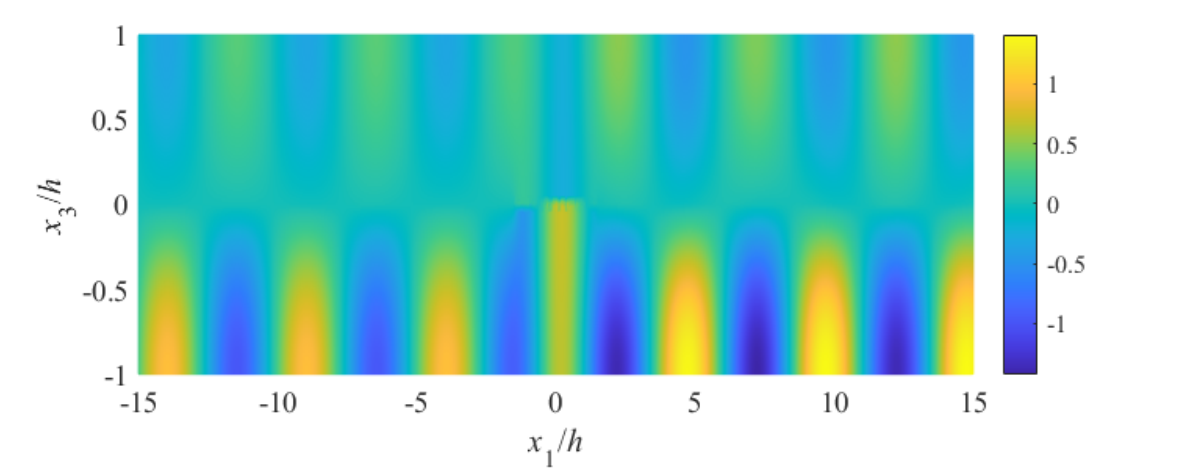}}
\caption{The total field displacement where the frequency is fixed at 2 MHz, the thickness of each layer is assumed as $0.5h$, and the second guided mode is selected as the incident wave.}
\label{fig17}
\end{figure}

In this section, various lengths of delamination ($h$, $1.2h$ and $1.5h$) are considered, where an aluminum-steel plate is adopted as the model, the thickness of each layer is assumed as $0.5h$, and the second guided mode is selected as the incident wave. It is clearly observed from Fig.\ref{fig16} that all curves have the same fluctuation when the new guided modes appear due to dispersion curves Fig.\ref{figmaterial_disper}, and the absolute value of reflection and transmission coefficients due to mode 2 is larger compared to mode 1 which indicates that the second guided mode accounts for the main component. The length of the delamination has a significant effect on the reflection coefficient and therefore the reflection coefficients of these guided modes can be used to detect. It should also be noted from Fig.\ref{fig17} that the total field displacement of different length of delamination has no noticeable difference.

\subsubsection{Different locations of the interface}
\begin{figure}[h]
	\centering
	\includegraphics[scale=0.55]{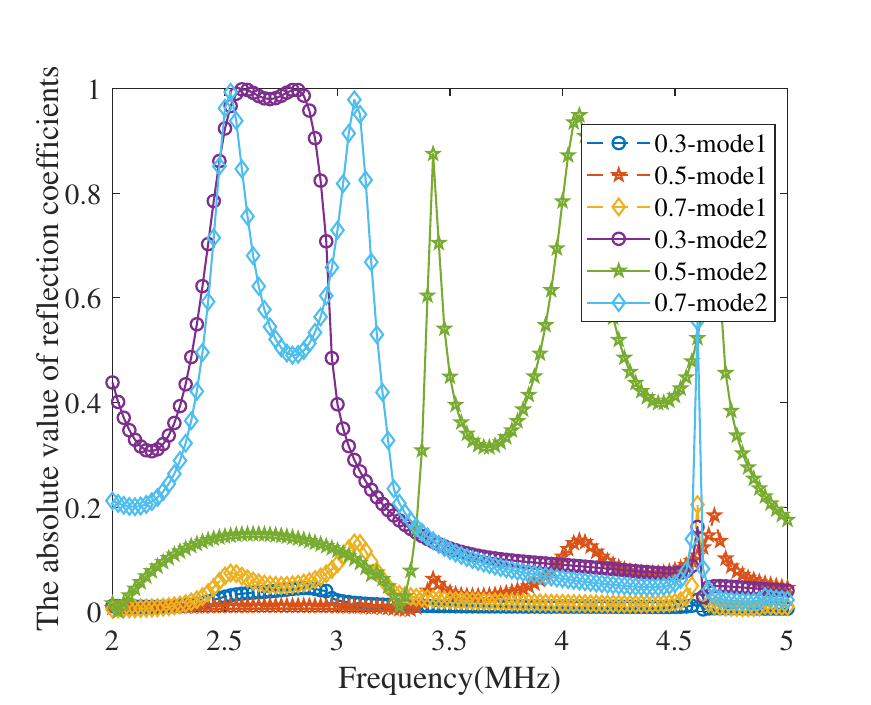}
    \includegraphics[scale=0.55]{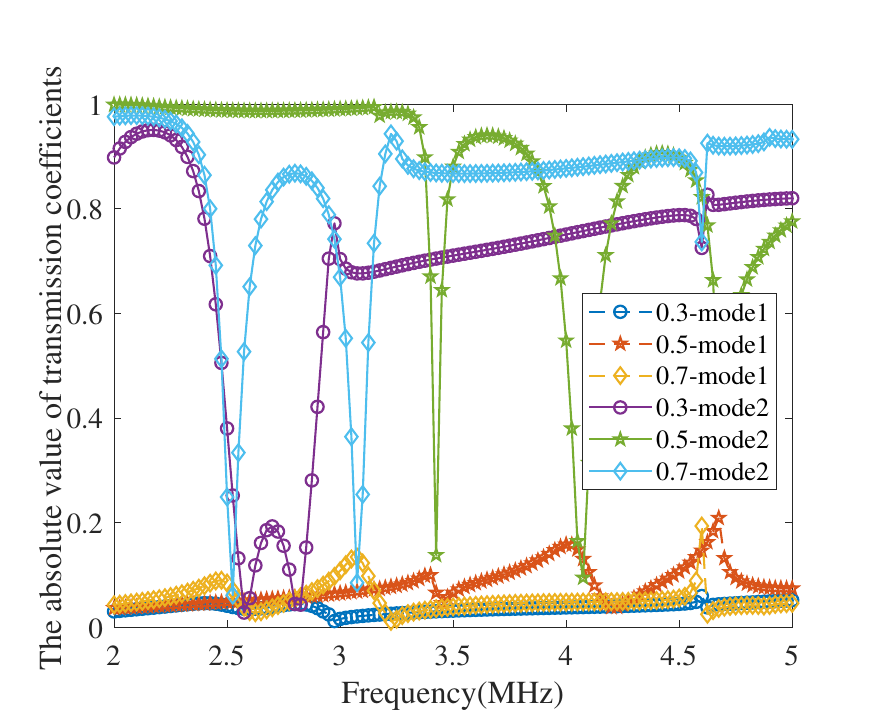}
	\caption{The absolute value of reflection and transmission coefficients varies with different location of the interface (different ratios of the first layer thickness and the total thickness (0.3, 0.5 and 0.7)) where an aluminum-steel plate is adopted as the model, the total thickness is assumed as $h$, and the second guided mode is selected as the incident wave}
	\label{fig18}
 \end{figure}

\begin{figure}[h]
\centering
 \subfigure[Ratio of the first layer thickness and the total thickness is 0.3;]{\includegraphics[scale=0.7]{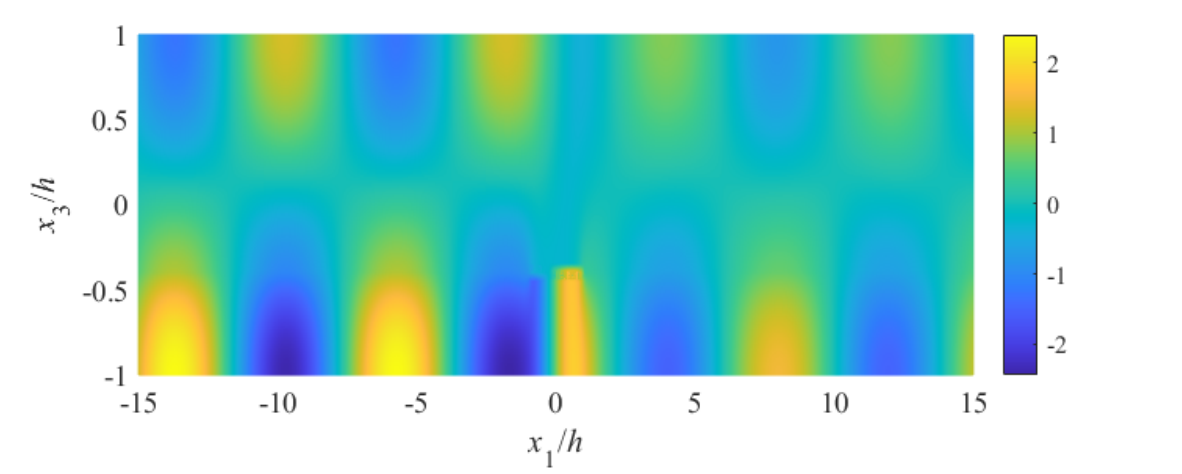}}
 \quad
\subfigure[Ratio of the first layer thickness and the total thickness is 0.7.]{\includegraphics[scale=0.7]{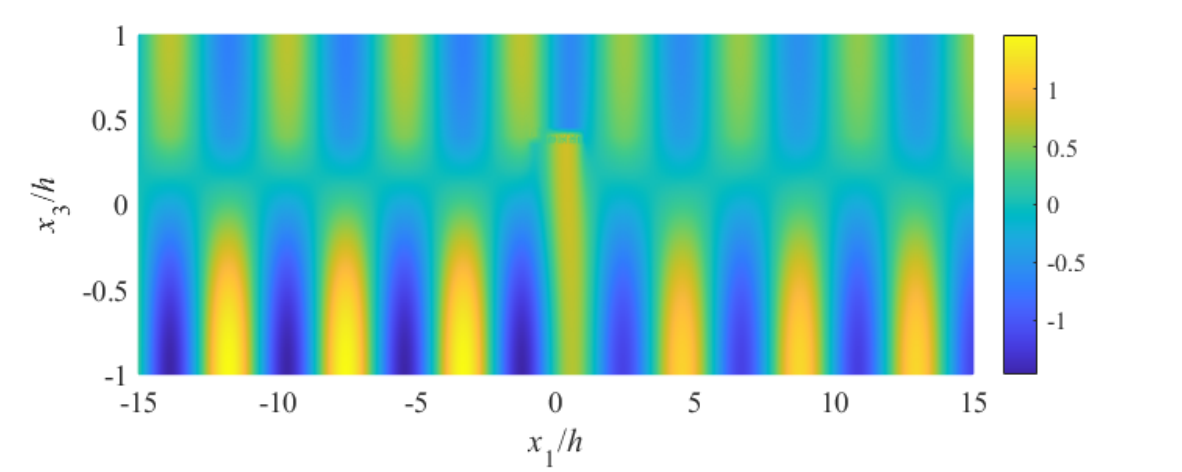}}
\caption{The total field displacement where the frequency is fixed at 2 MHz, the total thickness is assumed as $h$, and the second guided mode is selected as the incident wave.}
\label{fig19}
\end{figure}

In this section, different ratios of the first layer thickness and the total thickness (0.3, 0.5 and 0.7) are considered, where an aluminum-steel plate is adopted as the model, the total thickness is assumed as $h$, and the second guided mode is selected as the incident wave. It is clearly observed from Fig.\ref{fig18} that all curves have the same fluctuation when the new guided modes appear, and the absolute value of reflection and transmission coefficients due to mode 2 is larger compared to mode 1 which indicates that the second guided mode accounts for the main component. The location of the interface has a significant effect on the reflection coefficient where the reflection curves due to ratios 0.3 and 0.7 have a large fluctuation in the low frequency range while the reflection curve due to ratio 0.5 is relatively flat. It should also be noted from Fig.\ref{fig19} that the total field displacements of different ratios have noticeable difference.

\section{Conclusion}

This paper proposed a two-dimensional DtN FEM to analyze the scattering of SH guided waves due to an interface delamination in a bi-material plate. This method is simple and elegant, which has advantages on dimension reduction and needs no absorption medium or perfectly matched layer to suppress the reflected waves compared to traditional FEM. The conclusion can be specified:

1.	Normalized displacements at both ends converges to 1, which indicates that as the propagation distance increases, the non-propagating modes decay rapidly, leaving only the guided wave modes. Therefore, this scattering phenomenon coincides perfectly with the previously defined far-field assumption. The maximum error of energy balance is only around 0.25$\%$ and the average error is lower than $0.05\%$ which ensures sufficient calculation accuracy. The reflection and transmission coefficients computed by FEM are perfectly consistent with BEM's results which verifies the correctness and accuracy of our proposed DtN FEM.;

2.	In parametric analysis of different materials, the new mode of case Al-St is the first to appear and 
both the absolute values of reflection and transmission coefficients of mode 1 due to case Al-St are the largest compared with other two cases. The total displacement of case St-St is almost anti-symmetric about the interface since the anti-symmetric mode 2 accounts for the main component, and the geometry and material property are both symmetric about the interface;

3.	In parametric analysis of various lengths of the delamination, the length of the delamination has a significant effect on the reflection coefficient and therefore the reflection coefficients of these guided modes can be used to detect. The total field displacement of different length of delamination has no noticeable difference;

4.	In parametric analysis of different locations of the interface, the location of the interface has a significant effect on the reflection coefficient where the reflection curves due to ratios 0.3 and 0.7 have a large fluctuation in the low frequency range while the reflection curve due to ratio 0.5 is relatively flat. It should also be noted that the total field displacements of different ratios have noticeable difference.

For the future development, the scattering data from forward analysis by proposed DtN FEM will subsequently be used for the inverse analysis of reconstructing both the location and length of the delamination.

\section{Conflicts of Interest}

The authors declare that they have no known competing financial interests or personal relationships that could have appeared to influence the work reported in this paper.

\section*{Acknowledgments}

This work was supported by the China Scholarships Council (No. 202106830041). 


\bibliography{mybibfile}

\end{document}